\documentclass[a4paper, 12pt, reqno]{amsart}
\usepackage{amsmath, amssymb, caption, latexsym, graphicx}
\usepackage[ansinew]{inputenc}
\usepackage{fontenc}
\usepackage[dvips, dvipsnames, usenames]{color}
\usepackage[stable]{footmisc}
\usepackage[backref]{hyperref}
\hypersetup{colorlinks=true, urlcolor=blue, citecolor=blue}
\usepackage{fancyhdr}

\newcommand{\cnoindent}{\vspace{12pt}\noindent}
\newtheorem{proposition}{Proposition}

\newcommand{\s}[1]{$\langle #1_{i} \rangle_{i \in \mathbb{N}}$}

\newcommand{\w}{$\omega$ }

\newcommand{\ms}{\: }

\newcommand{\worder}{$\omega$-or\-der }
\newcommand{\wordered}{$\omega$-or\-de\-red }

\newcommand{\worderp}{$\omega$-or\-der}

\newcommand{\worderingp}{$\omega$-or\-de\-ring}

\begin{document}

\title{Hilbert's machine and \LARGE{\w}\normalsize{-order}}

\maketitle

\begin{center}

\small{
Antonio León\\
I.E.S. Francisco Salinas, Salamanca, Spain\\
}

\href{http://www.interciencia.es}{http://www.interciencia.es}\\
\href{mailto:aleon@interciencia.es}{aleon@interciencia.es}

\end{center}

\fancyhead{}

\pagestyle{fancy}

\fancyhead[CE, CO]{\small Hilbert's machine}

\renewcommand{\headrulewidth}{0.2pt}

\fancyfoot{}

\fancyfoot[LE, RO]{ \thepage}

\begin{abstract}
Hilbert's machine is a supertask machine inspired by Hilbert's Hotel whose functioning
leads to a contradictory result involving $\omega$-ordering and then the actual infinity.
\end{abstract}


\section{Introduction}

\noindent In the next discussion we will make use of a supermachine inspired by the
emblematic Hilbert's Hotel. But before beginning, let us relate some of the prodigious
(and suspicious) abilities of the illustrious Hotel. Its director, for instance, has
discovered a fantastic way of getting rich: it demands one euro to the guest of the room
1; this guest recover his euro by demanding one euro to the guest of the room 2; the
guest of the room 2 recover his euro by demanding one euro to the guest of the room 3,
and so on. Finally all guests recover his euro, and then our crafty director demands a
second euro to the guest of the room 1 which recover again his euro by demanding one euro
to the guest of the room 2, which recover again his euro by demanding one euro to the
guest of the room 3, and so on and on. Thousands of euros coming from the (infinitist)
nothingness to the pocket of our fortunate director.

\cnoindent Hilbert's Hotel is even capable of violating the laws of thermodynamics making
it possible the functioning of a perpetuum mobile: in fact we would only have to power
the appropriate machine with the calories obtained from the successive rooms of the
prodigious hotel in the same way its director got his euros. It is in fact unbelievable
that infinitists justify all those pathologies, and many other, in behalf of the
\emph{peculiarities} of the actual infinity. It is unbelievable that they prefer to
assume the pathological behaviour of the world before questioning the consistency of the
pathogene. But in the next discussion we will come to a contradiction that cannot be
easily subsumed in the picturesque nature of the actual infinity, a contradiction from
which it is impossible to escape.

\section{Hilbert's machine}

\noindent In the following conceptual discussion we will make use of a theoretical device
that will be referred to as \emph{Hilbert's machine}, composed of the following elements:

    \begin{enumerate}

        \item An infinite magnetic wire which is divided into two infinite parts, the
            left and the right side:

        \begin{enumerate}
            \item The right side is divided into an \wordered sequence of
                adjacent magnetic sections \s{s} which are indexed from left to
                right as $s_1$, $s_2$, $s_3$, $\dots$. They will be referred to
                as right sections.

            \item The left side is also divided into an \wordered sequence of
                adjacent sections \s{s'} indexed now from right to left as
                $\dots$, $s'_3$, $s'_2$, $s'_1$; being $s'_1$ adjacent to $s_1$.
                They will be referred to as left sections.

        \end{enumerate}

        \item An \wordered sequence of magnetic sliding beads \s{b} which are
            inserted in the magnetic wire as the beads of an abacus, being each bead
            $b_i$ initially placed on the right section $s_i$.

        \item A magnetic multidisplacement mechanism which moves simultaneously each
            bead exactly one section to the left, so that the bead placed on
            $s_{k,\ms k>1}$ is placed on $s_{k-1}$, the one placed on $s_1$ is placed
            on $s'_1$, and if one were placed on $s'_{k}$ it would be placed on
            $s'_{k+1}$. This simultaneous displacement of all beads \s{b} one section
            to the left will be termed \emph{magnetic multidisplacement}, or simply
            \emph{multidisplacement}. Multidisplacements are the \emph{only actions}
            Hilbert's machine can perform.

    \end{enumerate}

\noindent Let us now consider the following definition: we will say that a bead $b_i$ is
removed from the wire if, and only if, it is placed out of the wire as a consequence of a
multidisplacement. Although the impossibility of being removed from the infinite magnetic
wire\footnote{The last left magnetic section of the wire does not exist.} would
facilitate our discussion, we will assume that removal is, nevertheless, possible.
Consequently, we will impose the following restriction to the functioning of Hilbert's
machine: the machine will perform a multidisplacement if, and only if, the
multidisplacement does not remove any bead from the wire nor alters the original \worder
of the beads $b_1$, $b_2$, $b_3$, \dots (Hilbert's restriction).

\cnoindent Assume now that Hilbert's machine performs a magnetic multidisplacement $m_i$
at each one of the countably many instants $t_i$ of any \wordered sequence of instants
\s{t} defined within any finite interval of time $(t_a, t_b)$, for instance the classical
one defined by:
    \begin{equation}\label{eqn:ht definition of <ti>}
        t_i = t_a + (t_b - t_a)\sum_{k=1}^{i}\frac{1}{2^k}, \ \forall i \in \mathbb{N}
    \end{equation}
whose limit is $t_b$. In these conditions, at $t_b$ our machine will have completed the
\wordered sequence of multidisplacements \s{m}, i.e. a supertask. As is usual in
supertask theory\footnote{See \cite{Bernadete1964}, \cite{Grunbaum1967},
\cite{Grunbaum2001c}, \cite{Laraudogoitia1996}, \cite{Laraudogoitia1997},
\cite{Laraudogoitia2001}, \cite{Earman1998}, etc.} we will also assume that
multidisplacements are instantaneous. Although it is irrelevant for our conceptual
discussion, we could also assume that each multidisplacement lasts a finite amount of
time, for instance each $m_i$ could take a time $1/(2^{i+1})$. It seems appropriate at
this point to emphasize the conceptual nature of the discussion that follows. We are not
interested here in discussing the problems derived from the actual performance of
supertasks in our physical universe, as would be the case of the length of the wire or
the relativistic restrictions on the speed of the magnetic multidisplacements and the
like\footnote{\cite{Earman1993}, \cite{Earman1996}, \cite{Laraudogoitia1996},
\cite{Laraudogoitia1997}, \cite{Laraudogoitia1998}, \cite{Earman1998},
\cite{Laraudogoitia1999}, \cite{Norton1999}, \cite{Alper1999}, \cite{Alper2000},
\cite{Laraudogoitia2002}, \cite{Atkinson2006}, \cite{Atkinson2006b}, etc.}. We will
assume, therefore, that Hilbert's machine works in a conceptual universe in which no
physical restriction limits its functioning. Our only objective here is to examine the
consistency of \worderingp.

\section{Performing the supertask}

\noindent Consider the \wordered sequence of instants \s{t} defined according to
(\ref{eqn:ht definition of <ti>}), and a Hilbert's machine in the following initial
conditions:
    \begin{enumerate}
        \item At $t_a$ the machine is at rest.
        \item At $t_a$ each bead $b_i$ is on the right section $s_i$.
        \item At $t_a$ each left section $s'_i$ is empty.
    \end{enumerate}
\noindent Assume that, if Hilbert's restriction allows it, this machine performs exactly
one magnetic multidisplacement $m_i$ at each one of the countably many instants $t_i$ of
\s{t}, and only at them, being those successive multidisplacements the only performed
actions.

\cnoindent An \wordered sequence is one in which there exists a first element and each
element has an immediate successor. Consequently no last element exists. Thus, \wordered
sequences are both complete (as the actual infinity requires) and uncompletable (in the
sense that no last element completes them). The objective of the following discussion is
just to analyze the consequence of completing the uncompletable \wordered sequence of
multidisplacements \s{m}. We begin by proving the following basic proposition which is
directly derived from assuming the existence of \wordered sequences as complete
totalities.

\begin{proposition}\label{prp:hm all mi observe Hilbert restriction}
    \worder makes it possible that all multidisplacements $m_i$ of the \w-ordered
    sequence of multidisplacements \s{m} observe Hilbert's restriction.
\end{proposition}

\begin{proof}
It is evident the first multidisplacement $m_1$ observes Hilbert's restriction. In fact,
according to \worder each right section $s_{i, \, i>1}$ has an immediate predecessor to
the left and each left section has an immediate successor to the left. On the other hand,
$s_1$ has its own immediate predecessor to the left: $s'_1$. Thus, $b_1$ can be moved to
$s'_1$ and each $b_{i; i>1}$ to $s_{i-1}$. Consequently $m_1$ does not remove any bead
from the wire nor alter the original \worder of the beads since each bead $b_i$ remain
succeeded by its original immediate successor $b_{i+1}$ because all of them are
simultaneously moved one section the left. Assume the first n multidisplacements observe
Hilbert's restriction. Since each multidisplacement moves each bead exactly one section
to the left, after performing these first n multidisplacements, $b_1$ will have been
placed on $s'_n$, $b_{i,\,1 < i \leq n }$ on $s'_{n-i+1}$ and $b_{i,\, i>n}$ on
$s_{i-n}$. All these sections have an immediate predecessor (or successor) to the left so
that all beads can be simultaneously moved one section to the left without removing any
bead from the wire nor altering the initial \worder of the beads because each beads $b_i$
remains succeeded by its original immediate successor $b_{i+1}$ for the same reasons
above. In consequence, multidisplacement $m_{n+1}$ also observes Hilbert's restriction.
We have just proved that $m_1$ observes Hilbert's restriction, and that if the first n
multidisplacements observe Hilbert's restriction, then $m_{n+1}$ also observes Hilbert's
restriction. Therefore, every multidisplacement $m_i$ observes Hilbert's restriction.
\end{proof}

\noindent It is now possible to prove the following two contradictory results:

\begin{proposition}\label{prp:hm at tb supertask performed}
    At $t_b$ the \wordered sequence of multidisplacements \s{m} has been completed.
\end{proposition}

\begin{proof}
    According to Proposition \ref{prp:hm all mi observe Hilbert restriction} all
    multidisplacements \s{m} observe Hilbert restriction. Consequently all of them
    can be performed by Hilbert's machine without removing any bed from the wire nor altering
    the initial \worder of the magnetic beads. We will prove now that at $t_b$ all
    multidisplacements have been carried out. For this, consider the one to one
    correspondence $f$ between \s{t} and \s{m} defined by:
        \begin{equation}
            f(t_i) = m_i, \ \forall i \in \mathbb{N}
        \end{equation}
    \noindent Being $t_b$ the limit of the \wordered sequence \s{t}, and taking into account
    that each multidisplacement $m_i$ takes place at the precise instant $t_i$, the one to
    one correspondence $f$ together with the assumed completeness of all \wordered sequences,
    prove that at $t_b$ all multidisplacements \s{m} have been carried out. At $t_b$,
    therefore, the \wordered sequence of multidisplacements \s{m} has been completed.
\end{proof}

\begin{proposition}\label{prp:hm at tb supertask not performed}
    At $t_b$ the \wordered sequence of multidisplacements \s{m} has not been completed.
\end{proposition}

\begin{proof}
    According to Hilbert restriction and Proposition \ref{prp:hm all mi observe Hilbert
    restriction} no bead is removed from the wire and the initial \worder of \s{b} is
    preserved. Let therefore $b_n$ be any bead. Since it has not been removed from the wire,
    at $t_b$ it must of necessity be in on one of its magnetic sections. Assume it is on
    $s_k$. Since each multidisplacement moves $b_n$ a section to the left, only a finite
    number $n - k$ of multidisplacements will have been carried out to move $b_n$ from its
    initial section $s_n$ to $s_k$. Assume now $b_n$ is on $s'_h$ at $t_b$. In this case, and
    for the same reason above, only a finite number $n + h$ of multidisplacements will have
    been performed. We can therefore conclude that a $t_b$ only a finite number of
    multidisplacements can have been performed. So, at $t_b$ the \wordered sequence of
    multidisplacements \s{m} has not been completed.
\end{proof}

\section{Consequences}

\noindent We have just proved that at $t_b$ the \wordered sequence of multidisplacements
\s{m} has and has not been completed\footnote{Although we will not do it here, it is
possible to derive other contradictory results from the functioning of Hilbert's
machine}. Obviously, Hilbert's machine is a conceptual device whose theoretical existence
and functioning is only possible under the assumption of \worderp, that legitimates the
\emph{completeness} of the \wordered sequences \s{s}, \s{s'}, \s{b}, \s{m} and \s{t}.
Furthermore, the contradictory Propositions \ref{prp:hm at tb supertask performed} and
\ref{prp:hm at tb supertask not performed} are formal consequences of Proposition
\ref{prp:hm all mi observe Hilbert restriction}, which in turn is a formal consequence of
\worderp. It is, therefore, \worder the cause of the contradiction between Propositions
\ref{prp:hm at tb supertask performed} and \ref{prp:hm at tb supertask not performed}.

\cnoindent We will come to the same conclusion on the inconsistency of \worder by
comparing the functioning of the above infinite Hilbert's machine (symbolically
$H_\omega$) with the functioning of any finite Hilbert machine with a finite number $n$
of both right and left sections (symbolically $H_n$); being, as in the case of
$H_\omega$, a sequence of $n$ magnetic beads initially placed in the right side of the
wire, each bead $b_i$ on the section $s_i$. In effect, it is immediate to prove that,
according to Hilbert's restriction, $H_n$ can only perform $n$ multidisplacements because
the $(n+1)$-th multidisplacement would remove from the wire the bead $b_1$ initially
placed on the first right section $s_1$ and placed on the last left section $s'_n$ by
multidisplacement $m_n$. Thus $m_{n+1}$ does not observe Hilbert restriction and the
machine halts before performing $m_{n+1}$. $H_n$ halts with each left section $s'_i$
occupied by the bead $b_{n-i+1}$ and all right sections empty, and this is all. No
contradiction is derived from the functioning of $H_n$. Thus for any natural number $n$,
$H_n$ is consistent. Only infinite Hilbert's machine $H_\omega$ is inconsistent.
Consequently, and taking into account that \worder is the only difference between
$H_\omega$ and $H_{n, \: \forall \, n \in \mathbb{N}}$, only \worder can be the cause of
the inconsistency of $H_\omega$.

\cnoindent What the above contradiction proves, therefore, is not that a particular
supertask is inconsistent. What it proves is the inconsistency of \worder itself. Perhaps
we should not be surprised by this conclusion. After all, an \wordered sequence is one
which is both complete (as the actual infinity requires) and uncompletable (there is not
a last element that completes it). On the other hand, and as Cantor proved
\cite{Cantor1897}, \cite{Cantor1955}, \worder is an inevitable consequence of assuming
the existence of denumerable complete totalities. An existence axiomatically stated in
our days by the Axiom of Infinity, in all axiomatic set theories including ZFC and BNG
\cite{Suppes1972}, \cite{Stoll1979}. It is therefore that axiom the ultimate cause of the
contradiction between Propositions \ref{prp:hm at tb supertask performed} and \ref{prp:hm
at tb supertask not performed}.

\providecommand{\bysame}{\leavevmode\hbox to3em{\hrulefill}\thinspace}
\providecommand{\MR}{\relax\ifhmode\unskip\space\fi MR }
\providecommand{\MRhref}[2]{%
  \href{http://www.ams.org/mathscinet-getitem?mr=#1}{#2}
} \providecommand{\href}[2]{#2}

\end{document}